\documentclass[10pt,leqno]{article}
\usepackage{amssymb,amsbsy,amsmath,amsfonts,amssymb,amscd, mathrsfs}
\usepackage{color, graphicx} 

\usepackage[english]{babel}
\usepackage[T1]{fontenc}
\usepackage{indentfirst}

\makeatletter
\@addtoreset{equation}{section}
\makeatother

\newtheorem{statement}{}[section]
\newtheorem{theorem}[statement]{Theorem}
\newtheorem{lemma}[statement]{Lemma}

\newtheorem{claim}[statement]{Claim} 

\newcommand\C{\mathbb C}

\newcommand\D{\mathbb D}

\newcommand\e{{\rm e}}

\newcommand\ind{{\rm 1\kern-.30em I}}
\newcommand\qed{\hfill $\square$}

\let\phi=\varphi

\newcommand\converge{\mathop{\longrightarrow}\limits} 
\newcommand\tpg{|\hskip -0,5 pt \|} 
\newcommand\tpd{\| \hskip - 0,5 pt |}

\title{\bf Approximation numbers of composition operators on the Hardy space of the ball and of the polydisk}
\author{\it Daniel Li, Herv\'e Queff\'elec, Luis Rodr{\'\i}guez-Piazza\footnote{Supported by a Spanish research project MTM 2012-05622.}}

\date{\footnotesize \today}

\begin{document}

\maketitle

\noindent{\bf Abstract.} \emph{We give general estimates for the approximation numbers of composition operators on the Hardy space on the ball $B_d$ and 
the polydisk $\D^d$.}

\bigskip

\medskip


\section{Introduction} 

This work is an attempt to investigate approximation numbers of composition operators on the Hardy space $H^{2} (\Omega)$ where $\Omega$ is an open subset 
of $\C^d$, i.e. when we work with $d$ complex variables instead of one. In fact, we will essentially consider the two cases when $\Omega = B_d$ is the unit ball of 
$\C^d$ endowed with its usual hermitian norm $\Vert z \Vert = \big(\sum_{j = 1}^d |z_j|^2 \big)^{1/2}$ and  $\Omega = \D^d$ is the unit ball of $\C^d$ 
endowed with the sup-norm $\Vert z \Vert_\infty = \sup_{j = 1}^d |z_j|$, that is when $\Omega$ is the unit polydisk of $\C^d$. In order to treat these two 
cases jointly, we will work in the setting of bounded symmetric domains. 
\par\smallskip

An interesting feature is that the rate of decay of approximation numbers \emph{highly depends on $d$}, becoming slower and slower as $d$ increases, which 
might lead to think that no compact composition operators exist for truly infinite-dimensional symbols. We will see in the forthcoming paper 
\cite{LQR-polydisque infini} that this is not the case. 
\par\goodbreak

\section{Notations and background} 

A \emph{bounded symmetric domain} of $\C^d$ is an open convex and circled subset $\Omega$ of $\C^d$ such that for every point $a \in \Omega$, there is 
an involutive bi-holomorphic map $u \colon \Omega \to \Omega$ such that $a$ is an isolated fixed point of $\sigma$ (equivalently, $u (a) = a$ and 
$u ' (a) = - id$ (see \cite{VIGUE}, Proposition 3.1.1). \'E. Cartan showed that every bounded symmetric domain of $\C^d$ is homogeneous, i.e. the group 
of automorphisms of $\Omega$ acts transitively on $\Omega$: for every $a, b \in \Omega$, there is an automorphism $u$ of $\Omega$ such that 
$u (a) = b$ (see \cite{VIGUE}, p.~250). The unit ball $B_d$ and the polydisk $\D^d$ are examples of bounded symmetric domains. \par

The Shilov boundary $S_\Omega$ of such a domain $\Omega$ is the smallest closed set $F \subseteq \partial \Omega$ such that 
$\sup_{z \in \overline{\Omega}} | f (z) | = \sup_{z \in F} | f (z) |$ for every function $f$ holomorphic in some neighborhood of $\overline \Omega$. For 
example, the Shilov boundary of the bidisk is $S_{\D^2} = \{ (z_1, z_2) \in \C^2 \, ; \ |z_1| = |z_2| = 1\}$, whereas, its usual boundary $\partial \D^2$ is 
$\{ (z_1, z_2) \in \C^2 \, ; \ |z_1|, |z_2| \leq 1 \text{ and } |z_1| = 1 \text{ or } |z_2| = 1 \}$; for the unit ball $B_d$, the Shilov boundary is equal to the 
usual boundary ${\mathbb S}^{d - 1}$ (\cite{CLERC}, \S~4.1). Equivalently (see \cite{CLERC}, Theorem~4.2), $S_\Omega$ is the set of the extreme points 
of the convex set $\overline \Omega$. \par

If $\sigma$ is the unique probability measure on $S_\Omega$ invariant by the automorphisms $u$ of $\Omega$ such that $u (0) = 0$, the \emph{Hardy space} 
$H^2 (\Omega)$ is the space of all complex-valued holomorphic functions $f$ on $\Omega$ such that:
\begin{displaymath} 
\| f \|_{H^2 (\Omega)} := \bigg( \sup_{0 < r < 1} \int_{S_\Omega} | f (r \xi )|^2 \, d \sigma (\xi) \bigg)^{1/2} 
\end{displaymath} 
(see \cite{HAHN-MITCHELL}).  It is a Hilbert space (see \cite{HAHN-MITCHELL-TAMS}). \par \smallskip

A Schur map, associated with $\Omega$, will be a \emph{non-constant} analytic self-map of $\Omega$ into itself. It will be called \emph{truly $d$-dimensional} 
if the differential $\phi ' (a) \colon \C^d \to \C^d$ is an invertible linear map for at least one point $a \in \Omega$. Then, by the implicit function Theorem, 
$\phi (\Omega)$ has non-void interior. We say that the Schur map $\phi$ is a \emph{symbol} if it defines a \emph{bounded} composition operator 
$C_\phi \colon H^2 (\Omega) \to H^2 (\Omega)$ by $C_\phi (f) = f \circ \phi$. \par

Let us recall that if any Schur function generates a bounded composition operator on $H^2 (\D^d)$ when $d = 1$, this is no longer the case as soon as 
$d \geq 2$, as shown for example by the Schur map $\varphi (z_1, z_2) = (z_1, z_1)$.  Indeed, if say $d = 2$, taking $f (z) = (z_1 + z_2)^n$, we see that
\begin{displaymath} 
\Vert f \Vert_{2}^2 = \sum_{k = 0}^n \binom{n}{k}^2 = \binom{2 n}{n} \sim \frac{4^n}{\sqrt {\pi n}} \, \raise 1pt \hbox{,} 
\end{displaymath} 
while:
\begin{displaymath} 
\Vert C_{\phi}f \Vert_{2} = \Vert (2 z_1)^n \Vert_2 = 2^n \, .
\end{displaymath} 
The same phenomenon occurs on $H^2 (B_d)$ (\cite{MCCL}; see also \cite{CSW} and \cite{CW}).
\par\medskip

If $H$ is a Hilbert space and $T \colon H \to H$ is a bounded linear operator, the \emph{approximation numbers} of $T$ are defined, for $n \geq 1$ by:
\begin{equation} \label{approx numbers} 
a_n (T) = \inf_{\text{rank}\, R < n} \| T - R\| \,.
\end{equation} 
One has $\| T \| = a_1 (T) \geq a_2 (T) \geq \cdots \geq a_n (T) \geq a_{n + 1} (T) \geq \cdots$, and $T$ is compact if and only if 
$a_n (T) \converge_{n \to \infty} 0$.\par

The approximation numbers have (obviously) the following ideal property: for every bounded linear operators $S, U \colon H \to H$, one has:
\begin{displaymath} 
\qquad \qquad a_n ( STU) \leq \| S \| \, \| U \| \, a_n (T) \, , \qquad n = 1, 2 \ldots \, .
\end{displaymath} 
\par\smallskip

For an operator $T \colon H^{2} (\Omega) \to H^{2} (\Omega)$ with approximation numbers $a_{n} (T) = a_n$,  we will introduce the non-negative numbers 
$0 \leq \gamma_d^{-} (T) \leq \gamma_d^{+} (T) \leq \infty$ defined by:
\begin{equation} \label{gamma} 
\gamma_d^{-} (T) = \liminf_{n \to \infty} \frac{\log 1/a_n}{ n^{1/d}} \quad  \text{and} \quad
\gamma_d^{+} (T) = \limsup_{n \to \infty} \frac{\log 1/a_n}{ n^{1/d}} \, \cdot
\end{equation} 
The relevance of those parameters to the decay of approximation numbers is indicated by the following obvious facts, in which $0 < c \leq C < \infty$ denote 
constants independent of $n$:
\begin{align}   
\gamma_d^{-} (T) > 0 \  \quad & \Longleftrightarrow \quad a_n \leq C \, \e^{- c \, n^{1/d}} \, , \quad n = 1, 2, \ldots \label{rele1} \\ 
\gamma_d^{+} (T) < \infty \quad & \Longleftrightarrow \quad a_n \geq c \, \e^{- C n^{1/d}} \, , \quad n = 1, 2, \ldots \, . \label{rele2} 
\end{align}
So, the positivity of $\gamma_d^{-} (T)$ indicates that $a_n$ is ``small'' and the finiteness of  $\gamma_d^{+} (T)$ indicates that $a_n$ is ``big''.
\par\medskip

As usual, the notation $A \lesssim B$ means that there is a constant $c$ such that $A \leq c \, B$ and $A \approx B$ means that $A \lesssim B$ and $B \lesssim A$.

\section{Lower bound}

The next theorem shows that the approximation numbers of composition operators cannot be very small. We have already seen that in the one-dimensional case 
in \cite{LIQUEROD}. The important fact here is that this lower bound depend highly of the dimension.\par

\begin{theorem} \label{lower} 
Let $\Omega$ be a bounded symmetric domain of $\C^d$ and $\phi \colon \Omega \to \Omega$ be a truly $d$-dimensional Schur map inducing a compact 
composition operator $C_\phi \colon H^2 (\Omega) \to H^2 (\Omega)$. Then, for some constants $0 < c \leq C < \infty$, independent of $n$, we have: \par
\begin{displaymath} 
\qquad \qquad \quad a_{n} (C_\phi) \geq c \, \e^{- C n^{1/d}}, \qquad \forall n \geq 1, 
\end{displaymath} 
that is
\begin{displaymath} 
\gamma_d^{+} (C_\phi) < \infty \, .
\end{displaymath} 
\end{theorem}
\medskip

For proving that, we shall use the following results, the first of which is due to D. Clahane \cite{CLAHANE}, Theorem~2.1 (and to B.~MacCluer \cite{MCCL} 
in the particular case of the unit ball $B_d$). \par

\begin{theorem} [D. Clahane] \label{barbara} 
Let $\Omega$ be a bounded symmetric domain of $\C^d$ and $\phi \colon \Omega \to \Omega$ be a holomorphic map inducing a compact composition 
operator $C_\phi \colon H^2 (\Omega) \to H^2 (\Omega)$. Then $\phi$ has a unique fixed point $z_0 \in \Omega$ and the spectrum of $C_\phi$ consists of 
$0$, and all possible products of eigenvalues of the derivative $\phi ' (z_0)$.
\end{theorem}

When $\phi$ is truly $d$-dimensional, $0$ cannot be an eigenvalue of $C_\phi$ since if $f \circ \phi = 0$, then $f$ vanishes on $\phi (\Omega)$ which have a 
non-void interior, and hence $f \equiv 0$. Note that $1$ is an eigenvalue, by taking the product of zero eigenvalue of $\phi ' (z_0)$.\par\smallskip

In fact, in our case, we will not need the existence of $z_0$, for we will force $0$ to be a fixed point by a harmless change of the symbol $\phi$.
\smallskip

\begin{lemma} 
Let $H$ be a complex Hilbert space and $T \colon H \to H$ be a compact operator with eigenvalues $\lambda_1, \ldots, \lambda_n, \ldots$, written in 
non-increasing order and with singular values $a_n$, $n = 1 , 2, \ldots \,$. Then:
\begin{equation} \label{hermann} 
|\lambda_{2 n}|^2 \leq a_1 \, a_n \, .
\end{equation}
\end{lemma}

Indeed, it suffices to apply an immediate consequence of Weyl's inequalities, namely $|\lambda_n| \leq (a_1 \cdots a_n)^{1/n}$, with $n$ changed into $2 n$, 
and square to get 
\begin{displaymath} 
|\lambda_{2 n}|^2 \leq (a_1 \cdots a_{2 n})^{1/n} \leq (a_{1}^n \, a_{n}^n)^{1/n} = a_1\, a_n \, .
\end{displaymath} 
\par

\begin{lemma} \label{binome} 
Let $N_p$ be the number of multi-indices $\alpha = (\alpha_1, \ldots, \alpha_d)$ such that $|\alpha| = \alpha_1 + \cdots + \alpha_d \leq p$.  Then, 
as $p$ goes to infinity:
\begin{equation} \label{multi} 
N_p \sim \frac{p^d}{d!} \, \cdot
\end{equation}
\end{lemma}
\noindent {\bf Proof.} Let $n_k$ be the number of multi-indices $(\alpha_1, \ldots, \alpha_d, \alpha_{d + 1})$ such that 
$\alpha_1 + \cdots + \alpha_d + \alpha_{d + 1} = k$. We have (see \cite{LQ}, page~498), classically, for $|t| < 1$:
\begin{displaymath} 
\sum_{p = 0}^\infty n_p t^p = 
\bigg( \sum_{\alpha_1 = 0}^\infty t^{\alpha_1} \bigg) \cdots \bigg( \sum_{\alpha_{d + 1}}^\infty t^{\alpha_{d + 1}} \bigg) 
 =  \bigg( \sum_{k = 0}^\infty t^k \bigg)^{d + 1} =  \frac{1}{(1 - t)^{d + 1}}\, ; 
\end{displaymath} 
hence 
\begin{displaymath} 
n_p = \binom{d + p}{p} \, .
\end{displaymath} 
But $N_p = n_p$, and hence:
\begin{displaymath} 
N_p = \frac{(d + 1) \cdots (d + p)}{p!} = \frac{(d + p)!}{p! \, d!} \sim \frac{p^d}{d!} \, \raise 1pt \hbox{,}
\end{displaymath} 
by Stirling's formula for example. \qed

\goodbreak
\begin{claim}
We may assume that $\phi (0) = 0$ and $\phi ' (0)$ is invertible. 
\end{claim} 

\noindent{\bf Proof.} Since $\phi$ is truly $d$-dimensional, there exists $a \in \Omega$ such that $\phi ' (a)$ is invertible. Since $\Omega$ is homogeneous, 
there exist two automorphisms $\Phi_a$ and $\Phi_{\phi (a)}$ of $\Omega$ such that $\Phi_a (0) = a$ and $\Phi_{\phi (a)} [\phi (a)] = 0$. Set 
$\psi = \Phi_{\phi (a)} \circ \phi \circ \Phi_a$. Then $\psi (0) = 0$. Now, every analytic automorphism $\Phi$ of $\Omega$ induces a bounded composition 
operator on $H^2 (\Omega)$ and $C_\Phi^{- 1} = C_{\Phi^{- 1}}$ (\cite{CLAHANE}, Theorem~3.1); hence we can write 
$C_\psi = C_{\Phi_a} \circ C_\phi \circ C_{\Phi_{\phi (a)}}$ and it follows that $C_\psi$, as $C_\phi$, is compact. The ideal property of approximation 
numbers implies that, for $n = 1, 2, \ldots \,$, one has: 
\begin{displaymath} 
\big( \| C_{\Phi_a} \| \, \| C_{\Phi_{\phi (a)}} \| \big)^{- 1} \, a_n (C_\phi) \leq a_n (C_\psi) 
\leq \| C_{\Phi_{a}} \| \, \| C_{\Phi_{\phi (a)}} \| \, \, a_n (C_\phi) \, ,
\end{displaymath} 
so $\gamma_d^{-} (C_\psi) = \gamma_d^{-} (C_\phi)$. Moreover, using the chain rule, we see that $\psi ' (0)$ is invertible, since $\phi ' (a)$ is. \qed 
\par\bigskip

\noindent {\bf Proof of Theorem~\ref{lower}.} Let $\mu_1, \ldots, \mu_d$ be the eigenvalues of $\phi ' (0)$ and set 
$\min_{1\leq j\leq d} |\mu_j| = \e^{- A} > 0$. By Theorem~\ref{barbara}, the eigenvalues $\lambda_1, \ldots, \lambda_n, \ldots$ of $C_\varphi$ are the 
numbers $\mu_{1}^{\alpha_1} \cdots \mu_{d}^{\alpha_d}$ rearranged in non-increasing order. By definition, we have 
$\lambda_{N_p} = \prod_{j = 1}^d \mu_j^{\alpha_j}$ for some $d$-tuple $\alpha = (\alpha_1, \ldots, \alpha_d)$ such that $|\alpha| \leq p$.  Therefore, 
$|\lambda_{N_p}| \geq \e^{- A| \alpha|} \geq \e^{- A p}$. If $M_p = [N_p] / 2$ where $[ \, . \, ]$ stands for the integer part, equation \eqref{hermann} 
gives:
\begin{displaymath} 
\e^{- 2 A p} \leq |\lambda_{N_p}|^2 \leq |\lambda_{2 M_p}|^2 \leq a_1 \, a_{M_p}. 
\end{displaymath} 
Since $M_p \sim C_d \, p^d$ in view of Lemma~\ref{binome}, inverting this relation and using the monotonicity of the $a_n$'s clearly gives the claimed result.
\qed \par\medskip
\goodbreak

\section{An alternative approach for the polydisk and the unit ball}

The previous proof of Theorem~\ref{lower} is essentially a ``functional analysis'' one. It is interesting to give a proof using complex analysis tools instead of 
functional analysis ones. Moreover, this approach will be useful for the example in Section~\ref{section example}. \par 

In the general case, we are not be able to do that, and we only do it for the polydisk. The same approach works for the unit ball, by using 
results of B. Berndtsson in \cite{BER}. To save notation, we will give the proof in the case $d = 2$ but it clearly works in any dimension $d$. We will make use of 
the following theorem of P.~Beurling (\cite{GA} p.~285), in which the word \emph{interpolation sequence} refers to the space $H^{\infty}$ of bounded analytic 
functions on $\Omega$ ($\Omega = \D$ or $\D^2$), the interpolation constant $M_S$ of the  sequence $S = (s_j)$ being the smallest number $M$ such that, for 
any  sequence $(w_j)$ of data satisfying $\sup |w_j| \leq 1$, there exists $f \in H^{\infty} (\Omega)$ such that $f (s_j) = w_j$ and $\| f \|_\infty \leq M$.
\par \smallskip 
\goodbreak

\begin{theorem} [P.~Beurling] \label{fils} 
Let $(z_j)$ be an interpolating sequence in the unit disk $\D$, with interpolation constant $M$. Then, there exist analytic functions $f_{j}$, $j \geq 1$, on $\D$ 
such that:
\begin{displaymath} 
f_{j} (z_k) = \delta_{j, k} \qquad \text{and} \qquad \sum_{j = 1}^\infty |f_{j} (z)| \leq M \, , \quad  \forall z \in \D \, . 
\end{displaymath} 
As a consequence, if $A = (a_j)$ and $B = (b_k)$ are interpolation sequences of $\D$ with respective interpolation constants $M_A$ and $M_B$, their 
``cartesian product'' $(p_{j,k})_{j, k} = \big( (a_j, b_k) \big)_{j, k} $ is an interpolation sequence, with respect to $H^{\infty}(\D^2)$, with interpolation 
constant $\leq M_A \, M_B$.  
\end{theorem}
\goodbreak 

The consequence was observed in the paper \cite{BECHLI}. Indeed, if $(f_j)$ and $(g_k)$ are P.~Beurling's functions associated to $A$ and $B$ 
respectively, any sequence $(w_{j, k})$ with $\sup_{j, k}|w_{j, k}| \leq 1$ can be interpolated by the bounded analytic function 
\begin{displaymath} 
f (z, w) = \sum_{j, k\geq 1} w_{j, k}f_{j} (z) \, g_{k} (w)
\end{displaymath} 
which satisfies $\Vert f \Vert_\infty \leq M_A \, M_B$. \par\smallskip

Alternatively, in the sequel, we might use the result of \cite{BECHLI} on the sufficiency of  Carleson's condition on products of Gleason distances in the case of 
several variables. But we will stick to the previous approach. We now make use of the following lemma of \cite{LIQURO} which was enunciated in the 
one-dimensional case, but whose proof  works word for word in our new setting; indeed, the space of multipliers of $H^2 (\D^2)$ is (isometrically) 
$H^\infty (\D^2)$ and then one shows that the unconditionality constant of the sequence $(K_{s_j})_{1 \leq j \leq n}$ of reproducing kernels associated to 
a finite sequence $S = (s_j)_{1 \leq j \leq n}$ is less than $M_u$ (see also \cite{LIQUEROD}). Also note that he reproducing kernel of $H^2 (\D^2)$ is now, 
for $a = (a_1, a_2) \in \D^2$: 
\begin{displaymath} 
K_{a} (z_1, z_2) = \frac{1}{(1 - \overline{a_1} z_1) (1 - \overline{a_2} z_2)} \, \raise 1pt \hbox{,} 
\end{displaymath} 
with  $\Vert K_a \Vert^2 = [ (1 - |a_1|^2) (1 - |a_2|^2) ]^{- 1}$.  \par

\begin{lemma} \label{mu} 
Let $\phi \colon \D^2 \to \D^2$ be a symbol inducing a compact composition operator $C_\phi \colon H^2 (\D^2) \to H^2 (\D^2)$. Let 
$u = (u_1, \ldots, u_N)$ be a finite sequence of distinct points of $\D^2$ with interpolation constant $M_u$ and let $v_j = \varphi (u_j)$, $1 \leq j \leq N$. Let 
$M_v$ be the interpolation constant of $v = (v_1, \ldots, v_N)$. Then, setting:
\begin{displaymath} 
\mu_{N}^2 = \inf_{1 \leq j \leq N}\frac{| K_{v_j} |^2}{| K_{u_j} |^2} 
= \inf_{1 \leq j \leq N} \frac{(1 - |u_{j, 1}|^2) (1 - |u_{j, 2}|^2)}{(1 - |v_{j, 1}|^2) (1 - |v_{j, 2}|^2)} \, \raise 1pt \hbox{,}
\end{displaymath} 
with $u_j = (u_{j, 1}, u_{j, 2})$ and  $v_j = (v_{j, 1},  v_{j,2})$, one has: 
\begin{equation} \label{finland} 
a_{N} (C_\varphi) \geq c' \, \mu_N \, M_{u}^{- 1} \, M_{v}^{- 1} \geq c' \, \mu_N \, M_{v}^{- 2} \, .  
\end{equation}
\end{lemma}
\smallskip

The last inequality $M_u \leq M_v$ is proved as follows: let $\sup|w_j| \leq 1$ and choose $f \in H^\infty$ such that $f (v_j) = w_j$ and 
$\Vert f \Vert_\infty \leq M_v$; then $g = f \circ \varphi \in H^\infty$ and satisfies $\Vert g \Vert_\infty \leq M_v$ and $g (u_j) = f(v_j) = w_j$. \qed
\par \medskip

It remains to choose $u$ and $v$ and to estimate the parameters of the lemma. \par\smallskip 

As in the first proof, we may assume that $\phi (0) = 0$ and that the differential $\phi ' (0)$ is invertible. \par \smallskip

Since $\phi ' (0)$ is invertible, the set $\phi (\D^2)$  contains a closed polydisk of radius $0 < r < 1$ with center $0$.  We then take for $v$ the sequence 
$v_{j, k} = (r \omega^j, r \omega^k)$ where $\omega$ is a primitive $n$th-root of unity, e.g. $\omega = \e^{2 i \pi / n}$. We have $v = A \times A$ where 
$A = (r \omega, r \omega^2, \ldots, r \omega^n)$ so that the sequence $v$ has length $N = n^2$. We know (\cite{GA}, p.~284) that $M_A = r^{1 - n}$, so that 
Theorem~\ref{fils} gives us $M_v \leq r^{2 - 2 n}$. We now write $v_j = \varphi (u_j)$ with $|u_j| \leq r$, which is always possible by decreasing $r$ if 
necessary (this $r$ can be ridiculously small, but remains positive). Finally,
\begin{displaymath} 
\frac{\Vert K_{v_j} \Vert^2}{\Vert K_{u_j} \Vert^2} \geq (1 - |u_{j, 1}|^2) (1 - |u_{j, 2}|^2) \geq  (1 - r^2)^2 \, . 
\end{displaymath} 
Collecting all those estimates and using \eqref{finland}, we obtain:
\begin{displaymath} 
a_{n^2} (C_\varphi) \geq (1 - r^2)^2 \, r^{4 n - 4} \geq c \, r^{4 n} \, .
\end{displaymath} 
Interpolating an arbitrary integer $m$ between two consecutive squares, we clearly obtain Theorem~\ref{lower} for $\D^2$ (note that in dimension $d$ a factor 
$(1 - r^2)^d$ instead of $(1 - r^2)^2$ shows up). \qed 
\par\medskip 

\goodbreak 


\section{An upper bound} 

Though the result of this section is undoubtedly true in the general setting of bounded symmetric domains, we are not familiar enough with complex analysis in 
several variables to work it out. Therefore, we will assume in this section that: 
\begin{equation} \label{class} 
\qquad \qquad \Omega = B_{l_1} \times \cdots \times B_{l_N} \, , \quad \text{with } l_1 + \cdots + l_N = d
\end{equation}
is the product of $N$ unit balls. That covers the case of the unit ball of $\C^d$ ($N = 1$) and the case of the polydisk of $\C^d$ ($N = d$ and 
$l_1 = \cdots = l_N = 1$). To save notations, we will assume in the sequel with $N = 2$.\par\medskip

A point $z = (z_j)_{1 \leq j \leq d} \in \Omega$ is of the form $z = (u, v)$ with  $u = (u_j)_{1 \leq j \leq l_1}$, $v = (v_j)_{l_1 < j \leq d}$ and 
$\sum_{j = 1}^{l_1}|u_j|^2 < 1$, $\sum_{j = l_1 + 1}^d |v_j|^2 < 1$. We see that $\Omega$ is the unit ball of $\C^d$ equipped with the following norm:  
\begin{equation} \label{omnorm} 
\tpg z \tpd = \max \bigg[ \bigg(\sum_{j = 1}^{l_1} |u_j|^2 \bigg)^{1/2}, \bigg(\sum_{j = l_1 + 1}^d |v_j|^2 \bigg)^{1/2} \bigg] \, ,
\end{equation}
where $z = (u, v)$ with $u \in \C^{l_1}$ and $v \in \C^{l_2}$.  \par\smallskip\goodbreak

The Shilov boundary of $\Omega$ is $S_\Omega = S_{l_1} \times S_{l_2}$ and the normalized invariant measure on $S_\Omega$ is 
$\sigma = \sigma_{l_1} \otimes \sigma_{l_2}$ where $\sigma_{l_1}$ and $\sigma_{l_2}$ denote respectively the area measure on the hermitian spheres 
$S_{l_1}$ and $S_{l_2}$. \par\medskip

The following is in Rudin (\cite{RUD} p.~16). \goodbreak
\begin{lemma} \label{valeur} 
The monomials $e_\alpha$, with $e_\alpha (z) = z^\alpha$, form an orthogonal basis of $H^2 (\Omega)$. Moreover if $\alpha = (\beta, \gamma)$ with 
$\beta = (\alpha_1, \ldots, \alpha_{l_1})$ and $\gamma = (\alpha_{l_1 + 1}, \ldots, \alpha_d)$,  then writing $z = (u, v)$ we have:
\begin{displaymath} 
\| e_\alpha\|^2 = \int_{S_{l_1} \times S_{l_2}}|u^{\beta}|^2 \,|v^\gamma|^2 \, d\sigma_{l_1} (u) \, d\sigma_{l_2} (v) 
= \frac{(l_1 - 1)! \, \beta!}{(l_1 - 1 +|\beta|)!} \, \frac{(l_2 - 1)! \, \gamma!}{(l_2 - 1 +|\gamma|)!} \, \cdot
\end{displaymath} 
Therefore, if $f = \sum _{\alpha} c_\alpha \, e_\alpha \in H^2(\Omega)$, one has:
\begin{displaymath} 
\| f \|^2 = \sum_{\alpha} |c_\alpha|^2 \frac{(l_1 - 1)! \, \beta!}{(l_1 - 1 + |\beta|)!} \, \frac{(l_2 - 1)! \, \gamma!}{(l_2 - 1 + |\gamma|)!} \, \cdot
\end{displaymath} 
\end{lemma}

We can now state the main result of that section, in which we set $ \| \phi \|_\infty := \sup_{z \in \Omega} \tpg \phi (z) \tpd$. 

\begin{theorem} \label{above} 
Let $\Omega = B_{l_1} \times B_{l_2}$, $d = l_1 + l_2$, and $\phi \colon \Omega \to \Omega$ be a truly $d$-dimensional Schur map, inducing a compact 
composition operator $C_\phi \colon H^2 (\Omega) \to H^2 (\Omega)$. Then, if $ \| \phi \|_\infty < 1$, one has $\gamma_d^{-} (C_\phi) > 0$, that is 
there exist some constants $0 < c \leq C < \infty$, independent of $n$, such that:
\begin{equation} 
\qquad \qquad a_{n} (C_\phi) \leq C \, \e^{- c n^{1/d}} \, , \qquad n = 1, 2, \ldots \, .
\end{equation} 
\end{theorem}

\noindent{\bf Proof.} Let us set  $r = \| \phi \|_\infty < 1$. Let $f = \sum c_\alpha \, e_\alpha \in H^2 (\Omega)$ with 
\begin{equation} \label{let} 
c_\alpha = \widehat{f} (\alpha) \text{ and } \| f \|^2 = \sum_{\alpha}|c_\alpha|^2 \| e_\alpha \|^2 \leq 1 \, .
\end{equation} 
Then $C_\phi f = \sum c_\alpha \phi^\alpha$. \par\smallskip

We approximate $C_\phi$ by the $N_n$-rank operator $R$ defined by 
\begin{displaymath} 
R f = \sum_{|\alpha| \leq n} c_\alpha \phi^\alpha 
\end{displaymath} 
and we  set 
$g = C_\phi (f) - R (f)$ as well as $\alpha = (\beta, \gamma)$ and $z = (u, v)$. We begin with  observing that 
$\frac{(l_1 - 1 + p)!}{(l_1 - 1)! p!} \leq  (p + 1)^{l_1 - 1}$ and  $\frac{(l_2 - 1 + q)!}{(l_2 - 1)! q!} \leq  (q + 1)^{l_2 - 1}$. Since 
$|c_\alpha| \leq \| e_\alpha \|^{- 1}$, we get by Lemma~\ref{valeur} and the multinomial formula: 
\begin{equation} \label{multinom} 
\sum_{|\beta| = p} \frac{p!}{\beta!} \, | \phi^{\beta} (u)|^2 = \bigg( \sum_{j = 1}^{l_1} | \phi_j (u)|^2 \bigg)^p 
\end{equation}
and a similar formula with $|\gamma|= q$ that, setting $p+q=N$:
\begin{align*}
\sum_{\substack{|\beta|= p \\ |\gamma|= q}} \| e_\alpha \|^{- 2} |\phi^{\alpha} & (z)|^2 
=\sum_{\substack{|\beta| = p \\ |\gamma|= q}} 
\frac{(l_1 - 1 + p)!}{\beta! (l_1 - 1)!} \frac{(l_2 - 1 + q)!}{\gamma! (l_2 - 1)!} \, | \phi^{\beta} (u)|^2 \, |\phi^{\gamma} (v)|^2 \\ 
& \leq (p + 1)^{l_1 - 1} (q + 1)^{l_2 - 1} \bigg(\sum_{j = 1}^{l_1}| \phi_j (u)|^2 \bigg)^p \bigg( \sum_{j = l_1 + 1}^d |\phi_j (v)|^2 \bigg)^q \\ 
& \leq (p + 1)^{l_1 - 1} (q + 1)^{l_2 - 1} \, r^{2 p} \, r^{2 q} \leq (N + 1)^{l_1 + l_2 - 2} \, r^{2 N} \, .
\end{align*}
We thus have for $z \in \Omega$ the pointwise estimate (where we used \eqref{let} and the Cauchy-Schwarz inequality):
\begin{displaymath} 
|g(z)|^2 \leq \sum_{|\alpha|> n} \| e_\alpha \|^{- 2}|\phi^\alpha (z)|^2 
\leq \sum_{N > n} \sum_{p + q = N} (N + 1)^{d - 2} r^{2 N} \leq C_d \, n^d \, r^{2 n}
\end{displaymath}
for all $z\in \Omega$. This now implies $\| (C_\phi - R) f \|_{H^2} = \| g \|_{H^2} \leq  C'_d \, n^{d/2} \, r^{n}$. Hence:
\begin{displaymath} 
\| C_\phi - R \| \leq  C'_d \, n^{d/2} \, r^n \, .
\end{displaymath} 
Therefore:
\begin{displaymath} 
a_{N_n + 1}  \leq  C'_d \, n^{d/2} \, r^n \, .  
\end{displaymath} 
Since $N_n \sim n^d$, we get, with $r < \rho < 1$:
\begin{displaymath} 
a_{n^d} \lesssim \rho^n \, .  
\end{displaymath} 
We end  the proof by interpolation between two indices of the form $n^d$. \qed

\goodbreak


\section{An example} \label{section example} 

For $0 < \theta < 1$, the lens map $\lambda_\theta$ of parameter $\theta$ is defined by:
\begin{equation} 
\lambda_\theta (z) = \frac{(1 + z)^\theta - (1 - z)^\theta}{(1 + z)^\theta + (1 - z)^\theta} 
\end{equation} 
(see \cite{Shapiro-livre} or \cite{LLQR}). \par\medskip

Let $\lambda_1 = \lambda_{\theta_1}, \ldots, \lambda_d = \lambda_{\theta_d}$ be lens maps of parameters $0 < \theta_1, \ldots, \theta_d < 1$. We define 
a multi-lens map $\phi$ on the polydisk $\D^d$ as:
\begin{equation} \label{multi-lens} 
\phi (z_1, \ldots, z_d) = \big( \lambda_1 (z_1), \ldots, \lambda_d (z_d) \big) \, , 
\end{equation} 
for $(z_1, \ldots, z_d) \in \D^d$. We write it $\phi = \lambda_1 \otimes \cdots \otimes \lambda_d$. \par

Since we may replace $\theta_1, \ldots, \theta_d$ by $\max_k \theta_k$ or by  $\inf_k \theta_k$ without changing the results, we will assume in the sequel that 
$\theta_1 = \cdots = \theta_d = \theta$, and we will say that the multi-lens map $\phi = \phi_\theta$ has parameter $\theta$.

\begin{theorem} \label{lentille} 
Let $\phi$ be a multi-lens map with parameter $\theta$. Then, for positive constants $a, b, a', b'$ depending only on $\theta$ and $d$, one has:
\begin{equation} \label{bi} 
a' \, \e^ {- b' n^{1/(2 d)}} \leq a_{n}(C_\phi)\leq a \, \e^ {- b \, n^{1/(2 d + 1)}}
\end{equation}
In particular, $\gamma_d^{-} (C_\phi) = 0$ even though $C_\phi$ is all Schatten classes.
\end{theorem}

The exponent $1 / (2 d + 1)$ in the upper estimate should certainly be $1 / (2 d)$, but our method does not give it. \par\medskip 

\noindent {\bf Proof.} 1) Let us first show that $C_\phi$ is Hilbert-Schmidt (and hence compact). We know by \cite{Shapiro-livre}, \S~2.3, that each 
composition operator $C_{\lambda_k}$ is Hilbert-Schmidt. Since $(e_\alpha)_\alpha$ is an orthonormal basis of $H^2 (\D^d)$, one has: 
\begin{align*} 
\| C_\phi \|_{HS}^2 
& = \sum_\alpha \| C_\phi (e_\alpha) \|_{H^2 (D^d)}^2 = \sum_\alpha \| \phi^\alpha \|_{H^2 (D^d)}^2 \\
& = \sum_\alpha \|  \lambda_1^{\alpha_1} \otimes \cdots \otimes \lambda_d^{\alpha_d} \|_{H^2 (D^d)}^2 \\ 
& = \sum_\alpha \| \lambda_1^{\alpha_1} \|_{H^2 (\D)} \cdots \| \lambda_d^{\alpha_d} \|_{H^2 (\D)}^2 \, , \quad 
\text{by Fubini's Theorem} \\ 
& = \prod_{k = 1}^d \sum_{\alpha_k = 0}^\infty \| \lambda_k^{\alpha_k} \|_{H^2 (\D)}^2 
= \prod_{k = 1}^d \sum_{\alpha_k = 0}^\infty \| C_{\lambda_k} (e_{\alpha_k}) \|_{H^2 (\D)}^2 \\
& = \prod_{k = 1}^d \| C_{\lambda_k} \|_{HS}^2 < + \infty \, ;
\end{align*} 
hence $C_\phi$ is Hilbert-Schmidt. Since $\| C_{\lambda_k} \|_{HS} \leq \frac{K}{1 - \theta}$ for some constant $K$ (see \cite{LLQR}, Lemma~2.2), 
one gets:
\begin{displaymath} 
\| C_\phi \|_{HS} \leq \Big( \frac{K}{1 - \theta} \Big)^d \, . 
\end{displaymath} 
Since the approximation numbers are non-increasing, one has:
\begin{displaymath} 
n \, [a_n (C_\phi) ]^2 \leq \sum_{l = 1}^n [a_l (C_\phi]^2 \leq \sum_{l = 1}^\infty [a_l (C_\phi]^2 = \| C_\phi \|_{HS}^2 \, ;
\end{displaymath} 
hence:
\begin{equation} \label{estimation provisoire}
a_n (C_\phi) \lesssim \frac{1}{\sqrt{n} \, (1 - \theta)^d} \, \cdot
\end{equation} 

As in \cite{LLQR}, \S~2, this inequality improves itself, by the semi-group property of the lens maps: 
$\lambda_\theta \circ \lambda_{\theta'} = \lambda_{\theta \theta'}$. Indeed, multi-lens maps have the same property:
\begin{displaymath} 
\phi_\theta \circ \phi_{\theta'} = \phi_{\theta \theta'} \, ,
\end{displaymath} 
and hence, for $0 < \tau < 1$ and $ k = 1, 2, \ldots$:
\begin{displaymath} 
C_{\phi_\tau^k} = [C_{\phi_\tau}]^k \, .
\end{displaymath} 

Now, the approximation numbers satisfy the sub-multiplicative property: $a_{m + n - 1} (ST) \leq a_m (S) \, a_n (T)$. Since 
$a_{m + n} (ST) \leq a_{m + n - 1} (ST)$, this implies that $a_{k n} (T) \leq [a_n (T)]^k$ for $n, k \geq 1$. \par

For $k \geq 1$ to be choosen later, let $\tau = \theta^{1 / k}$. We get, using \eqref{estimation provisoire} with $\tau$ instead of $\theta$:
\begin{displaymath} 
a_{k n} (C_{\phi_\theta}) = a_{k n} (C_{\phi_\tau}^k ) \leq [a_n (C_{\phi_\tau}) ]^k \lesssim \bigg( \frac{1}{\sqrt{n}\, (1 - \tau)^d} \bigg)^k
\leq \bigg( \frac{k^d}{\sqrt{n}\, (1 - \theta)^d} \bigg)^k \, .
\end{displaymath} 
since $1 - \theta = 1 - \tau^k \leq k (1 - \tau)$. \par

Choosing now for $k$ the integer part of $\delta n^{1/ (2 d)}$, where $\delta > 0$ is small enough (namely $\delta < 1 - \theta$), we get that:
\begin{displaymath} 
a_{k n} (C_{\phi_\theta}) \lesssim \e^{- b_1 k} \lesssim \e^{ - b_2 n^{1 / (2 d)}} \, .
\end{displaymath} 
Changing notation, we fall on, for every $N \geq 1$:
\begin{displaymath} 
a_N (C_{\phi_\theta}) \lesssim \e^{- b \, N^{1 / (2 d + 1)}} \, . 
\end{displaymath} 

This implies that, for all $p > 0$, $\sum_{N = 1}^\infty [a_N (C_{\phi_\theta}) ]^p < \infty$ i.e. $C_{\phi_\theta}$ is in all Schatten classes $S_p$. \par
\medskip

2) To prove the lower bound, we will use Theorem~\ref{fils} and Lemma~\ref{mu}. \par

Let $\sigma > 0$ and, for $1 \leq j_k \leq N$, $1 \leq k \leq d$:
\begin{displaymath} 
u_{j_1, \ldots, j_d} = (1 - \e^{- j_1 \sigma}, \ldots, 1 - \e^{- j_d \sigma}) \, .
\end{displaymath} 
Let:
\begin{displaymath} 
v_{j_1, \ldots, j_d} = \phi (u_{j_1, \ldots, j_d}) = \big( \lambda_1 (1 - \e^{- j_1 \sigma}), \ldots, \lambda_d (1 - \e^{- j_d \sigma}) \big) \, .
\end{displaymath} 

By \eqref{finland}, one has, with $N = n^d$:
\begin{equation} \label{finland-bis} 
a_N (C_\phi) \geq c' \mu_N M_v^{- 2} \, .
\end{equation} 
Actually, if 
\begin{displaymath} 
\mu_{k, N} = \inf_{1 \leq j_k \leq N} \frac{1 - |1 - \e^{- j_k \sigma}|^2}{1 - |\lambda_k (1 - \e^{- j_k \sigma}) |^2} \, \raise 1 pt \hbox{,} 
\end{displaymath} 
one has:
\begin{displaymath} 
a_N (C_\phi) \geq c' \prod_{1 \leq k \leq d} \mu_{k, N} M_v^{- 2} \, .
\end{displaymath} 
On the other hand, if $M_{k, v}$ is the interpolation constant of the sequence 
\begin{displaymath} 
\big( \lambda_k (1 - \e^{-  \sigma}), \ldots, \lambda_k (1 - \e^{- N \sigma}) \big) \, , 
\end{displaymath} 
of points of $\D$, one has $M_v \leq M_{1, v} \cdots M_{d, v}$, by Theorem~\ref{fils}; hence:
\begin{displaymath} 
a_N (C_\phi) \geq c' \prod_{1 \leq k \leq d} \mu_{k, N} M_{k, v}^{- 2} \, .
\end{displaymath} 

But we proved in \cite{LQR-Hp} (see the proof of Proposition~2.6 there) that:
\begin{displaymath} 
\mu_{k, N} M_{k, v}^{- 2} \gtrsim \e^{ - \beta \sqrt{n}} 
\end{displaymath} 
for some constant $\beta > 0$ depending only on $\theta$. We get hence:
\begin{displaymath} 
a_N (C_\phi) \gtrsim \e^{ - \beta d \, \sqrt{n}} \, .
\end{displaymath} 
Since $N = n^d$, we get, by interpolation, that, for every $N \geq 1$:
\begin{displaymath} 
a_N (C_\phi) \gtrsim \e^{ - \beta d \, N^{1 / (2 d)}} \, ,
\end{displaymath} 
and that ends the proof of Theorem~\ref{lentille}. \qed

\goodbreak

\newpage

\noindent
{\rm Daniel Li}, Univ Lille Nord de France, \\
U-Artois, Laboratoire de Math\'ematiques de Lens EA~2462 \\ 
\& F\'ed\'eration CNRS Nord-Pas-de-Calais FR~2956, \\
Facult\'e des Sciences Jean Perrin, Rue Jean Souvraz, S.P.\kern 1mm 18, \\
F-62\kern 1mm 300 LENS, FRANCE \\ 
daniel.li@euler.univ-artois.fr
\medskip\goodbreak

\noindent
{\rm Herv\'e Queff\'elec}, Univ Lille Nord de France, \\
USTL, Laboratoire Paul Painlev\'e U.M.R. CNRS 8524 \& 
F\'ed\'eration CNRS Nord-Pas-de-Calais FR~2956, \\
F-59\kern 1mm 655 VILLENEUVE D'ASCQ Cedex, 
FRANCE \\ 
Herve.Queffelec@univ-lille1.fr
\smallskip

\noindent
{\rm Luis Rodr{\'\i}guez-Piazza}, Universidad de Sevilla, \\
Facultad de Matem\'aticas, Departamento de An\'alisis Matem\'atico \& IMUS,\\ 
Apartado de Correos 1160,\\
41\kern 1mm 080 SEVILLA, SPAIN \\ 
piazza@us.es\par


\begin{thebibliography}{99} 

\bibitem {BAFILIQU} F.~Bayart, C.~Finet, D.~Li and H.~Queff\'elec, \emph{Composition operators on the Wiener-Dirichlet algebra}, 
J. Operator Theory 60 (2008), no.~1, 45--70.

\bibitem{BER} B.~Berndtsson, \emph{Interpolating sequences for $H^\infty$ in the ball}, Nederl. Akad. Wetensch. Indag. Math. 47 (1985), no. 1, 1--10.

\bibitem {BECHLI} B.~Berndtsson, S.-Y.~Chang, and K.-C.~Lin, \emph{Interpolating sequences in the polydisc}, Trans. Amer. Math. Soc. 302 (1987), 
161--169.

\bibitem {CSW} J. A.~Cima, C. S.~Stanton, W. R.~Wogen, \emph{On boundedness of composition operators on $H^2( B_2)$}, 
Proc. Amer. Math. Soc. 91 (1984), no. 2, 217--222. 

\bibitem {CW}  J. A.~Cima, W. R.~Wogen, \emph{Unbounded composition operators on $H^2 (B_2)$}, 
Proc. Amer. Math. Soc. 99 (1987), no. 3, 477--483. 

\bibitem {CLAHANE} D.~D.~Clahane, \emph{Spectra of compact composition operators over bounded symmetric domains}, 
 Integral Equations Operator Theory 51 (2005), no.~1, 41--56.
 
\bibitem {CLERC} J.-L.~Clerc, \emph{Geometry of the Shilov boundary of a bounded symmetric domain}, J. Geom. Symmetry Phys. 13 (2009), 25--74.
 
\bibitem {COMA} C.~Cowen, B.~MacCluer, Composition Operators on Spaces of Analytic Functions, Studies in Advanced Mathematics,  CRC Press (1994). 

\bibitem {GA} J.~Garnett,  Bounded Analytic Functions, Revised First Edition, Springer (2007).

\bibitem {HAHN-MITCHELL-TAMS}  K. T.~Hahn, J.~Mitchell,  \emph{$H^p$ spaces on bounded symmetric domains}, Trans. Amer. Math. Soc. 146 (1969), 
521--531.

\bibitem {HAHN-MITCHELL}  K. T.~Hahn, J.~Mitchell,  \emph{$H^p$ spaces on bounded symmetric domains}, Ann. Polon. Math. 28 (1973), 89--95.

\bibitem {LLQR} P. Lef\`evre, D. Li, H. Queff\'elec, L. Rodr{\'\i}guez-Piazza, \emph{Some new properties of composition operators associated to lens maps}, 
Israel J. Math. 195 (2) (2013), 801--824.

\bibitem {LQ} D. Li, H. Queff\'elec, Introduction \`a l'\'etude des espaces de Banach. Ana\-lyse et probabilit\'es,  
Cours Sp\'ecialis\'es 12, Soci\'et\'e Math\'ematique de France, Paris (2004).

\bibitem {LIQUEROD} D.~Li, H.~Queff\'elec, L.~Rodr{\'\i}guez-Piazza, \emph{On approximation numbers of composition operators}, 
J. Approx. Theory 164 (2012), no.~4, 431--459.

\bibitem {LIQURO} D.~Li, H.~Queff\'elec, L.~Rodr{\'\i}guez-Piazza, \emph{Estimates for approximation numbers of some classes of composition operators 
on the Hardy space}, Ann. Acad. Sci. Fenn. Math. 38 (2013), no.~2, 547--564.

\bibitem {LQR-Hp} D.~Li, H.~Queff\'elec, L.~Rodr{\'\i}guez-Piazza, \emph{Approximation numbers of compo\-sition operators on $H^p$}, 
submitted. \\ 
\texttt{https://hal-univ-artois.archives-ouvertes.fr/hal-01119589}

\bibitem {LQR-polydisque infini} D.~Li, H.~Queff\'elec, L.~Rodr{\'\i}guez-Piazza, \emph{Composition operators on the Hardy space of the infinite polydisk}, 
{\it in preparation}.

\bibitem {MCCL} B.~MacCluer, \emph{Spectra of compact composition operators on $H^{p} (B_N)$}, Analysis 4 (1984), 87--103.

\bibitem {RUD} W.~Rudin, Function Theory in the unit ball of $\C^n$, Second Edition, Springer (2008).

\bibitem {Shapiro-livre} J. H. Shapiro, Composition operators and classical function theory, 
Universitext, Tracts in Mathematics, Springer-Verlag, New-York (1993). 

\bibitem {VIGUE} J.-P.~Vigu\'e, \emph{Le groupe des automorphismes analytiques d'un domaine born\'e d'un espace de Banach complexe. Application aux 
domaines born\'es sym\'etriques},  Ann. Sci. \'Ecole Norm. Sup. (4) 9 (1976), no. 2, 203--281. 

\end{thebibliography}
\end{document}